\documentclass[tran,10pt]{IEEEtran}

\usepackage{amsmath,amssymb,float,arydshln,color}
\usepackage{psfrag,setspace,wrapfig,subfigure}
\usepackage[latin1]{inputenc}
\usepackage{enumerate}
\usepackage{empheq}

\usepackage{dsfont}
\usepackage[lined,ruled,commentsnumbered]{algorithm2e}
\usepackage{epsfig}
\usepackage[outdir=./]{epstopdf}
\usepackage{graphicx}

\usepackage{hyperref}
\usepackage{amsfonts}
\usepackage{cite}
\usepackage{url}

\usepackage{hhline}
\usepackage[table]{xcolor}
\usepackage{multirow}
\usepackage{tabu}
\usepackage{booktabs}
\usepackage{threeparttable}
\allowdisplaybreaks

\def\tran{^{\mathsf{T}}}

\def\sgn{{\rm sgn}}

\DeclareMathOperator*{\argmin}{arg\,min}
\newcommand{\bp}{ \begin{proof}}
\newcommand{\ep}{\end{proof} }

\newcommand{\Ex}{\mathbb{E}\hspace{0.05cm}}

\newcommand{\bm}[1]{\mbox{\boldmath $#1$}}

\newcommand{\be}{\begin{equation}}
\newcommand{\ee}{\end{equation}}
\newcommand{\bal}{\begin{align}}
\newcommand{\eal}{\end{align}}
\newcommand{\bq}{\begin{eqnarray}}
\newcommand{\eq}{\end{eqnarray}}
\newcommand{\bqn}{\begin{eqnarray*}}
\newcommand{\eqn}{\end{eqnarray*}}
\newcommand{\nn}{\nonumber}
\newcommand{\ba}{\left[ \begin{array}}
\newcommand{\ea}{\\ \end{array} \right]}

\newcommand{\define}{\;\stackrel{\Delta}{=}\;}

\def\bgamma {{\boldsymbol \gamma}}

\def\h{{\boldsymbol{h}}}

\def\x{{\boldsymbol{x}}}

\def\real{{\mathbb{R}}}

\def\Zint{{\mathchoice{\setbox1=\hbox{\sf Z}\copy1\kern-.75\wd1\box1}
{\setbox1=\hbox{\sf Z}\copy1\kern-.75\wd1\box1}
{\setbox1=\hbox{\scriptsize\sf Z}\copy1\kern-.75\wd1\box1}
{\setbox1=\hbox{\scriptsize\sf Z}\copy1\kern-.75\wd1\box1}}}

\makeatletter
\def\hlinewd#1{%
  \noalign{\ifnum0=`}\fi\hrule \@height #1 \futurelet
   \reserved@a\@xhline}
\makeatother
\title{Online Dual Coordinate Ascent Learning}

\author{\IEEEauthorblockN{Bicheng Ying, Kun Yuan, and Ali H. Sayed\\}
	\vspace{0.4cm}
	\IEEEauthorblockA{\large Department of Electrical Engineering\\
	University of California, Los Angeles} \thanks{This work was supported in part by NSF grants CCF-1524250 and ECCS-1407712, and by DARPA project N66001-14-2-4029.  Emails:\{ybc,kunyuan,sayed\}@ucla.edu}
	}
\renewcommand{\footnoterule}{%
	\kern -1pt
	\hrule width 1.1in height 1.2pt
	\kern 2pt
}
	
\IEEEoverridecommandlockouts
\begin{document}
\def\helvetica{phvr7t.tfm}
\def\helveticaoblique{phvro7t.tfm}
\def\helveticabold{phvb7t.tfm}
\def\helveticaboldoblique{phvbo7t.tfm}

\font\sfb=\helveticabold
=\helveticaboldoblique
\maketitle

\begin{abstract}
The stochastic dual coordinate-ascent (S-DCA) technique is a useful alternative to the traditional stochastic gradient-descent algorithm for solving large-scale optimization problems due to its scalability to large data sets and strong theoretical guarantees. However, the available S-DCA formulation is limited to finite sample sizes and relies on performing multiple passes over the same data. This formulation is not well-suited for online implementations where data keep streaming in. In this work, we develop an {\em online} dual coordinate-ascent (O-DCA) algorithm that is able to respond to streaming data and does not need to revisit the past data. This feature embeds the resulting construction with continuous adaptation, learning, and tracking abilities, which are particularly attractive for online learning scenarios.
\end{abstract}

%

\section{Introduction and Related Work}
We consider minimizing a regularized stochastic convex risk function of the form:
\be
\min_w J(w)\define \Ex Q(w;\x)+\rho R(w) \label{eq.1}
\ee
where the expectation is over the distribution of the data represented by the boldface letter $\x$,  $w\in\real^{M}$ is an unknown parameter vector, and $\rho\geq 0$ is a scaling factor. Moreover, the loss $Q(w;\x)$ is a  convex function over $w$ and it may be non-differentiable. The term $R(w)$ is a strongly-convex regularization factor such as $\ell_2$ or elastic-net regularization. In learning applications, it is customary for the data to consist of a scalar variable $\bm{\gamma}$ and an $M-$dimensional feature vector, $\h$, i.e., $\x=\{\bm{\gamma},\h\}$. We assume that the loss function depends on the data in the following manner:
\be
Q(w;\x) \define Q(\h\tran w;\bgamma)
\ee
This problem formulation is typical of many scenarios including  leaky least-mean-squares \cite{sayed2008adaptive}, support-vector machines \cite{cortes1995support,hsieh2008dual}, regularized logistic regression \cite{hastie2009elements,sayed2014adaptation}, and others.

In practice, it is customary to replace problem (\ref{eq.1}) by the minimization of a regularized {\em empirical} risk that is based on a collection of $N$ data points, namely,
\be
w^o\define \argmin_w\;\; \frac{1}{N}\sum_{n=1}^N Q\left(h_{n}\tran w; \gamma(n) \right) + \rho R(w)\label{eq.tradition_er}
\ee	
Here, the data $\{\gamma(n),h_n\}$ represent realizations arising from the distribution driving the variables $\{\bm{\gamma},\h\}$ and $N$ is the size of the data sample. One traditional, yet powerful, approach to solving problem \eqref{eq.tradition_er} is to employ a stochastic (sub)gradient method (SGD, for short) \cite{robbins1951stochastic,zhang2004solving, sayed2008adaptive, shalev2011pegasos}. This method can be implemented both in empirical form (involving repeated passes over a finite data sample) or online form (in response to streaming data). Either way, when a constant step-size, $\mu$, is used to drive the iteration, the SGD method can be shown to converge exponentially fast to a small neighborhood around $w^o$ with a steady-state error variance that is on the order of $O(\mu)$. Therefore, sufficiently small step-sizes can ensure satisfactory steady-state performance albeit at the cost of a slower convergence rate \cite{polyak1987introduction,sayed2014adaptation,bertsekas1989parallel}.

An alternative approach is to solve problem \eqref{eq.tradition_er} in the dual domain. Instead of minimizing \eqref{eq.tradition_er} directly, one can maximize the dual cost function using a coordinate-ascent algorithm\cite{luo1992convergence}. The dual problem involves maximizing over $N$ dual variables. Since updating all $N$ dual variables at each iteration can be costly, the coordinate-ascent implementation updates one dual variable at a time.
There have been several recent investigations along these lines in the literature with encouraging results. For example, references \cite{hsieh2008dual,yu2011dual} observed that a dual coordinate-ascent (DCA, for short) method can outperform the SGD algorithm when applied to large-scale SVM. Later, a {\em stochastic} version of DCA (denoted by S-DCA) was examined in
\cite{shalev2013stochastic,shalev2014accelerated} for more general risk functions. Compared with DCA, at each iteration, the stochastic implementation picks one data sample randomly (not cyclically) and updates the corresponding coordinate. Reference \cite{shalev2013stochastic} showed that S-DCA converges exponentially to the exact minimizer $w^o$ by running repeated passes over the finite data sample, which is a notable advantage over SGD.


Despite the apparent advantages in terms of theoretical guarantees and experimental performance, the stochastic DCA implementation suffers from three drawbacks for online scenarios. First, the available S-DCA implementation needs to know beforehand the size of the training data, $N$, since this value is explicitly employed in the algorithm. When data streams in, the value of $N$ is constantly changing and, therefore, the S-DCA implementation will not be applicable. Second, S-DCA needs to perform multiple passes over the same finite data sample. This situation is problematic for online operation when new data keeps streaming in and it is not practical to keep examining past data. 
Third, the S-DCA algorithm assigns the same weight ($1/N$) to each data sample in the training set. This is not ideal for scenarios where the minimizer $w^o$ can drift with time since it deprives the algorithm of adaptation and tracking abilities. 

In summary, while stochastic (sub)gradient techniques are able to solve problems of the type \eqref{eq.tradition_er} in an online manner, the available DCA and S-DCA algorithms lose this important feature. Motivated by these considerations, we focus in this article on developing an {\em online} stochastic coordinate-ascent algorithm, denoted by the letters O-DCA. While it shares some useful features with S-DCA, the online version allows the sample size to increase and continuously adjusts the weights that are assigned to the samples. The experimental results in this work illustrate the superior performance of O-DCA over SGD in terms of convergence rate and accuracy. We comment on these results by explaining how O-DCA shares interesting and revealing connections  with stochastic gradient-descent and stochastic proximal gradient algorithms in the {\rm primal} domain. Specifically, we will show that under $\ell_2-$regularization, the proposed O-DCA algorithm is related to a stochastic proximal gradient implementation, which helps explain the observed superior performance of O-DCA over SGD.

\section{Problem and Algorithm Formulation}
\subsection{Dual problem}

We first replace the empirical problem (\ref{eq.tradition_er}) by a more general weighted formulation that is able to capture several scenarios of interest as special cases. Namely, we consider instead the following problem:
\be
	\min_w\quad \frac{1}{\Delta_N}\sum_{n=1}^N \delta_{n,N}Q\left(h_{n}\tran w; \gamma(n) \right) + \rho R(w) \label{eq.primal_problem}
\ee
where $\delta_{n,N}\geq0$ is a weighting scalar factor and $\Delta_N>0$ is a normalization scalar factor. Both factors depend on the number $N$ of data points. Different choices for these factors correspond to different useful situations\cite{sayed2008adaptive,angelosante2010online}:
\begin{itemize}
	\item[(C1)] ({\bf Infinite-length window}): This case corresponds to the choice $\delta_{n,N}=1$ and $\Delta_N=N$, which reduces to (\ref{eq.tradition_er}). For these choices, all data starting from the remote past are scaled similarly.
	\item[(C2)]({\bf Exponential-weighting window}): In this case, we set
	\be
	\delta_{n,N} = \beta^{N-n}, \quad \Delta_N = \sum_{n=1}^{N} \beta^{N-n} = \frac{1-\beta^N}{1-\beta}
	\ee
	for some forgetting factor $\beta \in (0,1)$. Usually, the value of $\beta$ is very close to one, so that recent data are weighted more heavily than data from the remote past.
	
	\item[(C3)]({\bf Finite-length sliding window}): In this case, we focus on the most recent $L$ data points by setting $\Delta_N = L$ (for the initial stages when $N\leq L$, we set $\Delta_N = N$) and 
	\be
	\delta_{n,N} = \left\{
	\begin{aligned}
		1,\hspace{0.4cm}& {\rm when}\ n> N-L\\
		0,\hspace{0.4cm}& {\rm when}\ n\leq N-L\\
	\end{aligned}\right.
	\ee
\end{itemize}

%

In order to examine the dual problem of (\ref{eq.primal_problem}), we first rewrite it in the following equivalent form involving a set of linear constraints in terms of scalar variables $\{z(n)\}$:
\bq
	\hspace{-0.3cm}&\displaystyle\min_{w,\{z(n)\}}& \frac{1}{\Delta_N}\sum_{n=1}^N \delta_{n,N}Q\left(z(n); \gamma(n) \right) + \rho R(w) \label{cost function}\\
	\hspace{-0.3cm}&{\rm s.t.}& \frac{\delta_{n,N}}{\Delta_N}z(n) =  \frac{\delta_{n,N}}{\Delta_N} h_{n}\tran w, \hspace{0.3cm} n = 1, 2, \ldots, N \label{cst}
\eq
To see the equivalence, observe that if $\delta_{n,N}$ happens to be zero for index $n$, then it does not matter whether a constraint exists at that point in time or not because the corresponding loss term will disappear from the sum in (\ref{cost function}). 
Next, we introduce the Lagrangian function\cite{boyd2004convex}:\vspace{-0.2cm}
\begin{eqnarray}
\mathcal{L}(w,z,\lambda)
\hspace{-0.3cm}&=&\hspace{-0.3cm} \frac{1}{\Delta_N}\sum_{n=1}^N \Big(\delta_{n,N} Q\Big(z(n);\gamma(n) \Big) + \delta_{n,N}\lambda(n) z(n)\Big)\nn \\
&& {}+\rho\left(R(w) - \frac{1}{\rho\Delta_N}\sum_{n=1}^N \delta_{n,N}\lambda(n) h_{n}\tran w\right)
\end{eqnarray}
\noindent where $\{\lambda(n)\}$ are scalar Lagrange multipliers. Observe that we have as many Lagrange multipliers as the number of data samples and, therefore, their number increases continuously in the case of streaming data (which is the situation we are interested in). 
Next, we introduce the conjugate functions \cite{boyd2004convex}:
\bq
	Q^{\star}(x ; \gamma) &\define& \textstyle{\sup_z} \left\{x z - Q(z;\gamma)  \right\}\\
	R^{\star}(x)   &\define& \textstyle{\sup_w}  \left\{x\tran w - R(w) \right\}
\eq
We can now express the dual function, denoted by ${\cal D}_{N}(\lambda)$, in terms of these conjugate function as follows:
\begin{align}
	\mathcal{D}_N(\lambda) 
	=&\; \frac{1}{\Delta_N}\sum_{n=1}^N \delta_{n,N} \inf_{z(n)}\Big(Q\left(z(n) ; \gamma(n) \right) + \lambda(n) z(n)\Big)\nn \\
	&\;\;\;+{\rho} \inf_w \left(R(w) - \frac{1}{\rho\Delta_N}\sum_{n=1}^N\delta_{n,N} \lambda(n) h_{n}\tran w\right)\nn\\
	=&\; -  \frac{1}{\Delta_N}\sum_{n=1}^N \delta_{n,N} Q^{\star}\big(-\lambda(n) ;\gamma(n)\big)\nn\\
	&\quad\quad-\rho R^{\star}\left( \frac{1}{\rho\Delta_N}\sum_{n=1}^N\delta_{n,N} \lambda(n) h_{n}\right)\label{eq.dual_old}
\end{align}
From the infimum operation on the regularization term, the primal variable $w_N$ has to satisfy the following first-order optimality condition: 
\be
	\frac{1}{\rho\Delta_N}\sum_{n=1}^N\delta_{n,N} \lambda(n) h_{n} \in \partial R(w_N) \hspace{1cm}\label{eq.primal_first_order}
\ee
where $\partial R(w)$ denotes the sub-differential of $R(\cdot)$ at point $w$.  Now it is known that if a function $F(\cdot)$ is convex and closed, then it holds that\cite{bertsekas2003convex}:
\be
	x\in\partial F(y)\Longleftrightarrow y\in\partial F^\star(x) \label{lemma.subgradient}
\ee
Applying this property to (\ref{eq.primal_first_order}) and recalling that $R(w)$ is assumed to be strongly-convex, which implies its conjugate function, $R^{\star}(x)$, is continuously differentiable\cite{bertsekas2003convex}, we find that the primal variable $w_N$ is given by the following expression in terms of the gradient vector of the conjugate regularization function:
\be
	w_N = \nabla_x R^{\star}\left(w_N'\right)\label{eq.primal_two}
\ee
where we introduced the intermediate variable:
\bq
w'_N \define \frac{1}{\rho\Delta_N}\sum_{n=1}^N \delta_{n,N}\lambda(n) h_{n} \label{eq.primal_dual_var}
\eq
In Table \ref{table.reg} we list some common choices for the regularization term, its conjugate function, and gradient vector.
\begin{table}[t]
	\centering
	\caption{Typical choices for the regularization term, its conjugate function and gradient vector.}\vspace{-0.2cm}
	\begin{small}
	\begin{tabular}{|c|c|c|c|}
		\hline
		\cellcolor[gray]{0.8}&\cellcolor[gray]{0.8}  $R(w)$ &\cellcolor[gray]{0.8} $R^\star(x)$ &\cellcolor[gray]{0.8} $\nabla_x R^{\star}(x)$ \\
		\hhline{|====|}
		a)&${\color{white} \displaystyle\frac{1}{1}}\frac{1}{2}\|w\|^2$ & $\frac{1}{2}\|x\|^2$ & $x$\\
		\hline
		b)&${\color{white} \displaystyle\frac{1}{1}}\delta\|w\|_1 +\frac{1}{2}\|w\|^2$ & $\frac{1}{2} \|\mathcal{T}_{\delta}(x)\|^2$ & $\mathcal{T}_{\delta}(x)$\\
		\hline
		c)&${\hspace{-0.55cm}\color{white} \frac{\frac{\frac{1}{1}}{1}}{\frac{1}{1}}} \sum_{i=1}^M w(i)\log w(i)\hspace{-0.1cm}$ & $\hspace{-0.1cm}\log\left(\sum_{i=1}^Me^{x(i)}\right)\hspace{-0.12cm}$ & $\displaystyle\frac{ e^x}{\sum_i e^{x(i)}}$\\
		\hline
	\end{tabular}\label{table.reg}
	\begin{tablenotes}
		\item\hspace{-0.5cm} a) $\ell_2$-regularization.
		
		\item\hspace{-0.5cm} b) Elastic-net regularization. The entry-wise soft-threshold operator is defined as $[\mathcal{T}_{\delta}(w)]_i\hspace{-0.15cm}\define \hspace{-0.15cm}(|w(i)|-\delta)_+ \sgn(w(i))$ for the $i-$th entry.
		
		\item\hspace{-0.5cm} c) Regularization based on KL divergence. Here, $w(i)$ represents the $i-$th entry of $w$ and vector $w$ belongs to the probability simplex.\vspace{-0.2cm}
	\end{tablenotes}
	\end{small}
\end{table}
\noindent 
In this way, the dual function from (\ref{eq.dual_old}) can be expressed as:
\be
	\hspace{-0.1cm}\mathcal{D}_N(\lambda) = -\frac{1}{\Delta_N} \sum_{n=1}^N \delta_{n,N} Q^{\star}\big(-\lambda(n) ; \gamma(n)\big)-{\rho} R^\star(w'_N) \label{eq.dual_problem}
\ee
\subsection{Recursive constructions}
We still need to determine the dual variables, $\{\lambda(n)\}$, which help identify the primal solution through (\ref{eq.primal_two}) and (\ref{eq.primal_dual_var}). Before showing how to carry out this calculation, we observe first that expression (\ref{eq.primal_dual_var}) for the intermediate variable allows us to motivate recursive constructions for the primal variable. Revisiting the three scenarios we considered before:

\begin{itemize}
	\item[(C1)] {(\bf Infinite-length window)}: In this case, we have
	\begin{align}
		w'_N
		=&\; \frac{N-1}{N}w'_{N-1} +\frac{1}{\rho N}\lambda(N) h_{N}
	\label{eq.22}\end{align}
	\item[(C2)]({\bf Exponential-weighting window}): In this case, we have
	\begin{align}
	w'_N 
	=&\; \frac{\beta-\beta^{N}}{1-\beta^{N}}w'_{N-1}  +\frac{1-\beta}{\rho (1-\beta^N)}\lambda(N) h_{N}\label{exp.window}
	\end{align}
	\item[(C3)] ({\bf Finite-length sliding window}): We only consider the $N\hspace{-0.08cm}>\hspace{-0.08cm}L$ situation; the case $N \hspace{-0.08cm}\leq \hspace{-0.08cm}L$ can be handled similarly. 
	\begin{small}
		\begin{align}
		w'_N 
		=& \; w'_{N-1} - \frac{1}{\rho L}\lambda(N-L)h_{N-L} + \frac{1}{\rho L}\lambda(N)h_N\label{eq.25}
		\end{align}
	\end{small}
\end{itemize}

In all cases (\ref{eq.22})--(\ref{eq.25}), we find that there is a mapping that transforms $w_{N-1}'$ into $w_{N}'$. We denote this mapping generically by
\be
w_N' = f_N(w_{N-1}') +\alpha(N) \lambda(N) h_N \label{eq.primal_update}
\ee
for some scalar $\alpha(N)$ and function $f_N(x)$ given by (the function $f_N(\cdot)$ is affine in these three examples):
\be
f_N(x)=\left\{\begin{array}{ll} \frac{N-1}{N} x,&\hspace{-0.1cm}\mbox{\rm (infinite-length window)}\\
\frac{\beta-\beta^{N}}{1-\beta^N} x,&\hspace{-0.1cm}\mbox{\rm (exponential window)}\\
x-\frac{1}{\rho L}\lambda(N-L)h_{N-L},&\hspace{-0.1cm}\mbox{\rm (sliding window)}
\end{array}
\right.\nn
\ee
and
\be
\alpha(N) =\left\{\begin{array}{ll} \frac{1}{\rho N} ,&\mbox{\rm (infinite-length window)}\\
\frac{1-\beta}{\rho(1-\beta^N)},&\mbox{\rm (exponential window)}\\
\frac{1}{\rho L},&\mbox{\rm (sliding window)}
\end{array}
\right.
\ee
Observe that in the three cases considered above it holds that
\be		\alpha(N) =\frac{1}{\rho\Delta_N}
\ee

\subsection{Online algorithm}
Observe from (\ref{eq.dual_problem}) that the dual function ${\cal D}_N(\lambda)$ depends on both $N$ and $\lambda$, which creates a challenge for the development of an online algorithm. This is because the form of the dual function changes with $N$. Also, the number of dual variables increases with $N$. We therefore need an efficient method to seek the maximizer of the dual function. The main idea is as follows. When a new data point $\{\gamma(N),h_N\}$ streams in, we shall fix the previous Lagrange multipliers $\{\lambda(1),\lambda(2),\ldots,\lambda(N-1)\}$ at their existing values and then maximize ${\cal D}_N(\lambda)$ only with respect to $\lambda(N)$. It is important to emphasize that the motivation for this argument is somewhat different from traditional coordinate-ascent implementations. This is because the number of dual variables is now changing with time and, therefore, it is not possible to simply start from the solution of the last iteration. Instead, we extend the last solution into an enlarged vector that is one dimension higher and fix the leading entries of this longer vector to the dual variables from the last iteration. In this way,  we can write the dual function (\ref{eq.dual_problem}) as
\begin{align}
	\mathcal{D}_{N}(\lambda) \stackrel{\small (\ref{eq.primal_update})}{=}& - \frac{1}{\Delta_{N}} Q^{\star}\Big(-\lambda(N) ; \gamma(N)\Big) + \mbox{const}\nn\\
	 &\hspace{0.3cm}-\rho R^{\star}
	\left(f_N(w'_{N-1})+\alpha(N) \lambda(N) h_{N}  \right)
\end{align}
where the term ``const'' aggregates terms that are independent of $\lambda(N)$. By maximizing over $\lambda(N)$ we arrive at the proposed online dual coordinate-ascent (O-DCA) algorithm:
\begin{subequations}
\begin{eqnarray}
	\lambda(N) &=&\; \argmin_{\tau} \frac{1}{\Delta_{N}} Q^{\star}\Big(-\tau;\gamma(N)\Big)\label{eq.mainAlgorithm1}\\
	&&\hspace{1cm}{}+\rho R^{\star}
	\left( f_N(w'_{N-1})+ \tau\cdot\alpha(N)h_{N} \right)\nn \\
	w'_N &=&\;f_N(w'_{N-1})+\alpha(N)\lambda(N) h_{N}  \label{eq.mainAlgorithm2}\\
	w_N &=&\; \nabla R^{\star} (w'_N)\label{eq.mainAlgorithm3}
\end{eqnarray}
\end{subequations}
Observe that the algorithm involves three steps at each iteration $N$, when a new data $\{\gamma(N),h_N\}$ streams in. First, the optimal $\lambda(N)$ is determined by solving (\ref{eq.mainAlgorithm1}) Then, the intermediate estimate $w_N'$ is determined, followed by the evaluation of $w_N$. In comparison with the stochastic DCA (S-DCA) implementation of \cite{shalev2013stochastic,shalev2014accelerated}, three main differences stand out. First, at each iteration $N$, the proposed algorithm (\ref{eq.mainAlgorithm1})--(\ref{eq.mainAlgorithm3}) is employing a time-varying normalization factor $\Delta_N$, rather than a fixed $N$. This feature is critical for handling streaming data and to enable adaptation and tracking. Second, each data $\{\gamma(N),h_N\}$ is only used once, which is necessary for streaming data scenarios; multiple passes over the data are not practical in this case. And, third, more weighting is assigned to recent data than past data, which is important for scenarios with drifting minimizers. 
 
 In cases when the loss function $Q(\cdot)$ is non-differentiable, it is often helpful to smooth the output of O-DCA as follows:
\be
	\bar{w}_N  \define \frac{1}{S_N}\sum_{n=1}^N \kappa^{N-n} w_{n},\quad {\rm where}\;\; S_N \hspace{-0.1cm}\define \hspace{-0.1cm} \sum_{n=1}^{N}\kappa^{N-n}
\ee
and the weight factor $\kappa \hspace{-0.1cm}\in \hspace{-0.1cm} [0,1]$. Computing $\bar{w}_N$ can be implemented efficiently, e.g., by using the same recursive method used before in (\ref{exp.window}) to find that.
\subsection{Relation to Stochastic Primal Algrithms}\label{sectio.connection}
The online DCA algorithm (\ref{eq.mainAlgorithm1})--(\ref{eq.mainAlgorithm3}) that we just derived has strong connections with learning algorithms in the primal domain, especially when $\ell_2-$regularization is employed, i.e., $R(w)={1\over 2}\|w\|^2$. In this case, solving the argmin problem in (\ref{eq.mainAlgorithm1}) requires that we determine a $\lambda(N)$ that satisfies the following first-order condition:
\be
	 h_N\tran \nabla R_x^\star \Big(f_N(w'_{N-1})+\alpha(N)\lambda(N) h_{N} \Big) \in \partial Q^\star \Big(-\lambda(N);\gamma(N)\Big)
\ee
Using property (\ref{lemma.subgradient}) and update step (\ref{eq.mainAlgorithm2}), we conclude that:
\be
	\lambda(N) \in -\partial Q\Big(h_N\tran \nabla_x R^\star(w_N');\gamma(N)\Big)
\ee
It follows that we can rewrite the update for the intermediate variable in the O-DCA algorithm (\ref{eq.mainAlgorithm1})--(\ref{eq.mainAlgorithm2}) as follows:
\be
	\hspace{-0.1cm}w_N' = f_N(w_{N-1}') \hspace{-0.05cm}-\hspace{-0.05cm}
	{\alpha(N) \partial Q\left(h_N\tran \nabla R_x^\star(w_N') ; \gamma(N)\right) h_N}
\ee
For $\ell_2-$regularization we have $\nabla_{x} R^{\star}(x)=x$ and therefore the above recursion, along with the last equality (\ref{eq.mainAlgorithm3}), show that the proposed O-DCA algorithm reduces to the following insightful form (notice that $w_N$ appears inside the sub-differential instead of $w_{N-1}$):
\be
w_N^{\rm ODCA} = {f_N(w_{N-1}^{\rm ODCA})}-\alpha(N) \partial Q\Big(h_N\tran  w_N^{\rm ODCA};\gamma(N)\Big) h_N \label{eq.ODCA}
\ee
This form reveals important connections with iterative algorithms in the primal domain. Indeed, note that if we were to use an online stochastic (sub-)gradient descent (SGD) method to solve (\ref{eq.1}) with $\ell_2-$ regularization, we would have obtained the following recursion:
\be
w_{N}^{\rm SGD} = (1  -\mu \rho )w_{N-1}^{\rm SGD}
- \mu \partial Q\Big(h_N\tran w_{N-1}^{\rm SGD};\gamma(N) \Big) h_N\label{eq.SGD}
\ee
On the other hand, if we were to use an online stochastic proximal gradient (SPG) method 
 to solve (\ref{eq.1}), again with $\ell_2-$ regularization, we would have obtained:
\begin{align}
	w^{\rm SPG}_N
	=&\; {\rm prox}_{\mu Q} \Big((1-\mu \rho) w^{\rm SPG}_{N-1} \Big) \nn \\
	=&\; (1-\mu\rho)w^{\rm SPG}_{N-1} - \mu \partial Q\Big(h_N\tran w_{N}^{\rm SPG} ;\gamma(N)\Big) h_N \label{eq.SPG}
\end{align}
where in the last step we used the fact that $u = {\rm prox}_Q(x) \Longleftrightarrow x-u \in \partial Q(u)$. Comparing (\ref{eq.ODCA}), (\ref{eq.SGD}), and (\ref{eq.SPG}) under $\ell_2-$regularization, we observe that although O-DCA was formulated in the dual domain, it can still be viewed as one form of a proximal implementation with the variable $w_{N}^{\rm OCDA}$ appearing on the right-hand side of (\ref{eq.ODCA}) inside the sub-differential term, as happens with $w^{\rm SPG}_N$ in (\ref{eq.SPG}).

\section{Specific Loss functions}
We illustrate the above connections more explicitly by considering two important cases: least-mean-squares error designs and support vector machines. In both cases, for simplicity, we continue to employ $\ell_2-$regularization, $R(w)={1\over 2}\|w\|^2$.

\subsection{Least-Mean-Squares Learning}
In this case we have
\bq
	Q\Big(z(N) ; \gamma(N)\Big)\hspace{-0.2cm} &=& \hspace{-0.2cm} \frac{1}{2}\big(\gamma(N) - z(N)\big)^2\\
	Q^\star \Big(\lambda(N);\gamma(N)\Big)\hspace{-0.2cm} &=&\hspace{-0.2cm} \frac{1}{2}\lambda^2(N) + \gamma(N)\lambda(N)
\eq
Therefore, assuming an infinite-length window, we need to solve the following optimization problem  to find $\lambda(N)$:
\be
	\hspace{-0.1cm}\min_{\tau}\frac{1}{N}\left(\frac{1}{2}\tau^2 - \gamma(N)\tau\right) + \frac{\rho}{2} \left\|\frac{N-1}{N} w_{N-1}+\frac{\tau h_{N}}{\rho N}  \right\|^2
\ee
Setting the derivative relative to $\tau$ equal to zero at $\tau=\lambda(N)$ we find that 
the O-DCA algorithm (\ref{eq.mainAlgorithm1})--(\ref{eq.mainAlgorithm3}) reduces to:
\begin{subequations}
\begin{align}
	\hspace{-0.2cm}	\lambda(N) =&\;\left(1+\frac{\|h_{N}\|^2}{\rho N} \right) ^{-1} \Big(\gamma(N)-\frac{N-1}{N} h_{N}\tran w_{N-1} \Big)\\ 
	\hspace{-0.2cm}	w_{N} =&\;\frac{N-1}{ N }  w_{N-1}+\frac{\lambda(N) h_{N}}{\rho N}
\end{align}	\end{subequations}
Assuming $N$ is large enough,  we can merge the two equations into:
\begin{align}
	\hspace{-0.25cm}w_{N} =  w_{N-1}\hspace{-0.05cm}+\hspace{-0.05cm}\frac{1}{\rho N }h_{N}\Big(\gamma(N)\hspace{-0.05cm}-\hspace{-0.05cm} h_{N}\tran w_{N-1} \Big)\label{eq.41}
\end{align}
The O-DCA implementation (\ref{eq.41}) approaches a leaky-LMS implementation with a decaying step-size of the form $\mu_N=1/\rho N$, namely,
\be
	w_{N}  = (1-\rho\mu_N)w_{N-1} + \mu_N h_N(\gamma_N - h_N\tran w_{N-1})
\ee
\subsection{Support-Vector Machines}
In this case, we have (where $\gamma(N)=\pm 1$):
\be
	Q\Big(z(N) ; \gamma(N) \Big)= \max\{0,1-\gamma(N) z(N)\} \label{eq.hinge}
\ee
with conjugate function (in compact form):
\be
	\hspace{-0.3cm}Q^\star\big(\lambda(N);\gamma(N)\big) \hspace{-0.1cm}=\hspace{-0.1cm}
	\left\{
	\begin{aligned}
		\hspace{-0.1cm}\gamma(N)\lambda(N),\;&  \;\;\,{\rm if}\;\gamma(N)\lambda(N)\hspace{-0.1cm}\in\hspace{-0.1cm}[\hspace{-0.03cm}-\hspace{-0.03cm}1,0] \\
		\hspace{-0.1cm}+\infty,\hspace{1cm}& \;\;\,{\rm otherwise}
	\end{aligned}
	 \right.
\ee
In this example, we consider the exponential weighting window so that the dual variable is found by solving:
\begin{align}
	\lambda(N) \hspace{-0.1cm}=&\hspace{-0.1cm}\; \argmin_{\tau} \left\{-\frac{\gamma(N)\tau}{\Delta_{N}} + \frac{\rho}{2}\left\|\frac{\beta -\beta^N}{1-\beta^N} w_{N-1}+\frac{\tau h_{N}}{\rho\Delta_{N}} \right\|^2\right\}\nn\\
	&\; {\rm subject\ to}\quad \gamma(N)\tau\in[0,1] 
\end{align}
This optimization problem involves a truncated parabola function as a cost objective. The minimizer of the quadratic cost occurs at
\be
	\lambda(N) =\; \frac{\rho\Delta_N}{\|h_{N}\|^2}\left( \gamma(N)-\frac{\beta -\beta^N}{1-\beta^N}  h\tran_{N}w_{N-1}\right)
\ee
We still need to adjust this value, by means of a projection operation, in order to meet the constraint. To simplify the projection, we multiply both sides of the above relation by $\gamma(N)$ and use $\gamma^2(N)=1$ to find that the O-DCA algorithm (\ref{eq.mainAlgorithm1})--(\ref{eq.mainAlgorithm3}) reduces to:\vspace{-0.2cm}
\be
	\gamma(N)\lambda(N) = \frac{\rho\Delta_N}{\|h_{N}\|^2}\left( 1-\frac{\beta -\beta^N}{1-\beta^N}\gamma(N)  h\tran_{N}w_{N-1}\right)\hspace{1.7cm}
\ee\be
	w_{N} = \frac{\beta -\beta^N}{1-\beta^N}w_{N-1}+\frac{1-\beta}{\rho (1-\beta^N)} \gamma(N)h_{N}{\textstyle\prod_{[0,1]}} \big[\gamma(N)\lambda(N)\big]\nn\vspace{-0.1cm}
\ee
where $\Pi_{[0,1]}(a)$ projects the real number $a$ into the interval $[0,1]$. If we let $\beta = 1-\mu\rho$ and assume a small enough $\mu$ so that the value of $\beta$ is close to one, the above two equations can be merged into the following format when $N$ is large enough:
\begin{align}
	w_{N}^{\rm ODCA} =\; (1-\mu\rho) w_{N-1}^{\rm ODCA} + {} \label{eq.ODCA_SVM}\hspace{3.3cm}
\end{align}\vspace{-0.6cm}
\begin{align}
	\;\;\hspace{1cm} \mu\gamma(N) h_N \cdot{\textstyle\prod_{[0,1]}}\left[\frac{1}{\mu\|h_N\|^2}\left(1-\gamma(N)h_N\tran w^{\rm ODCA}_{N-1}\right)\right]\nn
\end{align}
For comparison purposes, we list the stochastic subgradient solution for SVM here\cite{ying2015performance,shalev2011pegasos}:
\be
	w_N^{\rm SGD} = (1-\mu\rho) w_{N-1}^{\rm SGD} + \mu\gamma(N) h_N\cdot\mathbb{I}[1-\gamma(N)h_N\tran w_{N-1}^{\rm SGD} \geq 0] \label{eq.SGD_SVM}
\ee
where $\mathbb{I}[\cdot]$ is the indicator function, which is equal to one when the argument is true and zero otherwise.  Comparing (\ref{eq.ODCA_SVM}) with (\ref{eq.SGD_SVM}), it is clear that O-DCA replaces the 
indicator function (which involves a sudden jump from $0$ to $1$) by a smoothed linear transition from $0$ to $1$ with slope proportional to $1/\mu$.

\section{Simulations}
We illustrate the performance of O-DCA using the hinge loss (\ref{eq.hinge}) applied to two datasets. The test data is obtained from the LIBSVM website\footnote{\hspace{-0.05cm}\url{http://www.csie.ntu.edu.tw/~cjlin/libsvmtools/datasets}}. We first use the Adult dataset with 11,220 training data and 21,341 testing data in 123 feature dimensions. In the figure, we compare our algorithm with the S-DCA method from \cite{shalev2013stochastic}, the stochastic sub-gradient method from \cite{ying2015performance}, and with LIBSVM\cite{CC01a}. To perform a fair comparison, we compare the performance based upon the number of iterations (except for LIBSVM). The parameter setting is as follows. All algorithms use $\rho = 0.001$. For O-DCA, we choose $\beta=0.99995$ and the step-size for the sub-gradient implementation is 0.05, so that it satisfies $\beta = 1-\mu\rho$. Since S-DCA is not designed for online learning, we feed S-DCA with one-fifth of the training data, and run over 5 epochs so that the  total number of iterations will match.
\begin{figure}[h]
	\vspace{-0.2cm}\epsfxsize 7cm \epsfclipon
	\begin{center}
		\leavevmode \epsffile{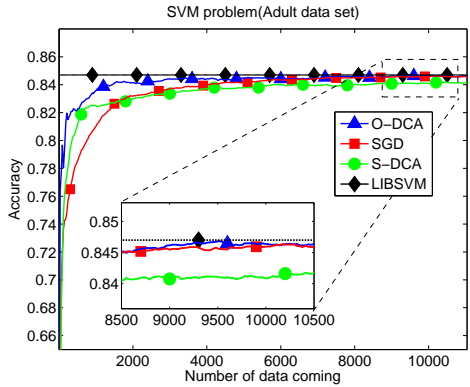}\vspace{-0.2cm}\caption{{\small Comparison of the performance of O-DCA, S-DCA, SGD, and LIBSVM on the Adult dataset in terms of iterations.}} \label{fig-adult.label}
	\end{center}
\end{figure}

The second dataset is the Reuters Corpus Volume I (RCV1) data with 20242 training data and 253843 testing data consisting of 47236 feature dimensions. Similarly, we set
$\rho=10^{-4}$, $\beta = 0.99998$ for O-DCA,  $\mu=0.2$ for sub-gradient so that $\beta = 1- \mu\rho$, and epochs~$=5$ for S-DCA.
\begin{figure}[h]
	\vspace{-0.32cm}
	\epsfxsize 7cm \epsfclipon
	\begin{center}
		\leavevmode \epsffile{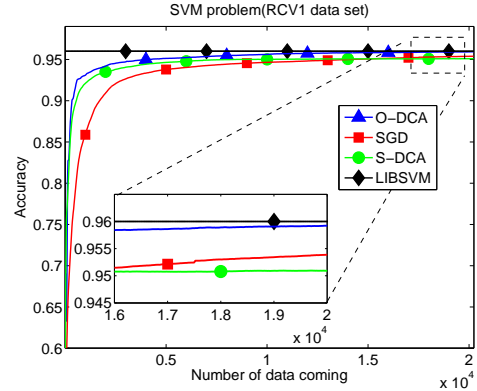}\vspace{-0.25cm}\caption{{\small Comparison of the performance of O-DCA, S-DCA, SGD, and LIBSVM on
	 the RCV dataset in terms of iterations.}} \label{fig-rcv.label}
	\end{center}
\end{figure}

\vspace{-0.3cm}
\bibliographystyle{IEEEbib}
\bibliography{ref_ODCA}

\end{document}